\documentclass[a4paper,12pt,leqno]{article}
\usepackage{amsmath,amsfonts,amscd,amsthm,ProtocoloEulerIH,url,txfonts,mathrsfs}
\usepackage[all]{xy}
\usepackage[colorlinks]{hyperref}
\usepackage{pdfsync}

\title{\Large\bf  {Equivariant intersection cohomology of the circle actions}\footnote{This work has been partially supported by the UPV/EHU grant EHU09/04 and by the Spanish MICINN grant MTM2010-15471.}}

\author{
Jos\'e Ignacio Royo Prieto\thanks{Departamento de Matem\'atica
Aplicada
UPV/EHU.
Alameda de Urquijo s/n.
48013 BILBAO,
SPAIN. \hspace{5cm}  {\sl joseignacio.royo@ehu.es}.}
\\
\\
{\small Universidad del Pa{\'\i}s Vasco/Euskal Herriko
Unibertsitatea} \and
Martintxo Saralegi-Aranguren
\thanks{Univ Lille Nord de France F-59\kern 1mm 000 LILLE, FRANCE.
UArtois, Laboratoire de Math\'ematiques de Lens EA~2462.
F\'ed\'eration CNRS Nord-Pas-de-Calais FR~2956.
Facult\'e des Sciences Jean Perrin.
Rue Jean Souvraz, S.P.\kern 1mm 18.
 F-62\kern 1mm 300 LENS, FRANCE.
{\sl saralegi@euler.univ-artois.fr}. }
\\
\\ {\small Universit\'e d'Artois }
}

\begin{document}

\maketitle

\noindent{\small Accepted for publication in
{\sc Revista de la Real Academia de Ciencias Exactas, F{\'\i}sicas y Naturales. Serie A. Matem\'aticas}, 2012.}

\noindent{\small 
DOI: 10.1007/s13398-012-0097-z. The final publication is available at {\tt springerlink.com}}.

\medskip

\noindent{\em To Heisuke Hironaka on the occasion of his 80th birthday.}

\begin{abstract}
Circle actions on pseudomanifolds have been studied  in \cite{PS} by using intersection cohomology (see also \cite{HS}). In this paper, we continue that study using a more powerful tool,
the equivariant intersection cohomology
\cite{Bry,J}.

In this paper, we prove  that the orbit space $B$ and the Euler class of the action $\Phi \colon \sbat \times X \to X$ determine both the equivariant intersection cohomology of the pseudomanifold $X$ and its localization.

We also  construct a spectral
sequence converging to the equivariant intersection cohomology of $X$ whose third
term is described in terms of the intersection cohomology of $B$.

 \end{abstract}

We consider an action $\Phi \colon \sbat \times X \to X$ of the circle on a pseudomanifold $X$ whose orbit space $B$ is again a pseudomanifold (cf. (1.1)).
We have seen in \cite{PS} that  the intersection cohomology
of $X$ is determined by  $B$ and  the Euler class $e \in \lau{\IH}{2}{\per{e}}{B}$. In this paper we prove that those two data determine some other structures. The main results of this work are the following:

\bigskip

$  \rightsquigarrow$
The equivariant intersection cohomology\footnote{See \cite{Bry,J}.} $\Hiru{\mbox{\boldmath $\IH$}}{\sbat}{X}$  of $X$  has a  $ \Lambda \eut$-perverse  algebra structure\footnote{
$ \Lambda \eut =\hiru{H}{*}{\C\P^\infty}$.}.
We prove that this structure is determined by
$B$ and the Euler class
$e \in \lau{\IH}{2}{\per{e}}{B}$ (cf. Proposition \ref{aritz}).

\bigskip
\addtocounter{footnote}{-2}
$  \rightsquigarrow$  The localization\footnotemark\  $\Hiru{\mbox{\boldmath $\IL$}}{\sbat}{X}$  of
$\Hiru{\mbox{\boldmath $\IH$}}{\sbat}{X}$ has a
perverse superalgebra structure. We prove that this structure is determined by $B$ and the Euler class
$e \in \lau{\IH}{2}{\per{e}}{B}$
(cf. Proposition \ref{bi}).
\addtocounter{footnote}{2}

\bigskip

$  \rightsquigarrow$ For each perversity $\per{p}$ we construct a spectral sequence converging  to $\lau{\IH}{*}{\per{p},\sbat}{X}$
whose third term is described in terms of
$B$ (cf. Proposition \ref{spec})\footnote{As C. Allday pointed out to us,
this spectral sequence degenerates into the
Skjelbred exact sequence of \cite{Sk} when $\per{p}=\per{0}$ (cf. Proposition \ref{skj}).}.

\medskip

In the last section, we illustrate the  results of this work with some particular examples.
In the sequel, any manifold will be considered  connected, second
countable, Haussdorff,
without boundary and smooth
(of class $C^\infty$).

The authors wish to thank the referee for the indications given in order to improve this paper.

\section{Modelled actions}

We recall in this section some fundamental notions about the objects we deal with. Namely, modelled actions on unfolded pseudomanifolds, the intersection cohomology and the Gysin sequence.

\prg {\bf Unfolded pseudomanifolds  (\cite{Pa}, \cite{S2})}

Let $X$ be a Hausdorff, paracompact and 2nd countable topological space.
We say that $X$ is a {\em stratified space} if it is provided with a stratification $\mathcal{S}$, that is, a finite partition by connected sets called {\em strata}, which satisfy the following condition: for any two strata $S,S'\in\mathcal{S}$, then $S\cap \overline{S'}\ne\emptyset$ implies $S\subset\overline{S'}$ (in this case, we will write $S\le S'$).

Notice that $(\mathcal{S},\le)$ is a partially ordered set, and that a stratum is maximal if and only if it is open. Such a stratum will be said to be {\em regular}. Every other strata is said to be {\em singular}. The union $\Sigma$ of every singular strata is called the {\em singular part} of the stratified space $X$. Its complement $X\backslash\Sigma$ is called the {\em regular part}.

Stratified pseudomanifolds were introduced by Goresky and MacPherson in order to extend the Poincar\'e duality to stratified spaces. A stratified space $X$ is a {\em stratified pseudomanifold} if for every stratum $S$ of $X$ there exists a family of {\em charts}, i.e., embeddings  $\alpha\colon U_{\alpha}\times c(L_S)\to X$
such that:
\begin{itemize}
\item $c(L_S)$ is the cone of a compact stratified space $L_S$ called the {\em link} of $S$;

\item $\{ U_{\alpha} \} _{\alpha}$ is an open cover of $S$;

\item $\alpha(u,\ast)=u$ for every $u\in U_{\alpha}$, being $\ast$ the apex of the cone $c(L_S)$.
\end{itemize}

Now, we recall the notion of an unfolding, which we use in order to define the intersection cohomology of a stratified pseudomanifold by means of differential forms. We will say that a stratified pseudomanifold $X$ is an {\em unfolded pseudomanifold} if it admits an {\em unfolding}, which consists of a manifold $\widetilde{X}$, a surjective, proper continuous function $\mathscr{L}\colon\widetilde{X}\to X$ and a family of unfoldings $\mathscr{L}_L\colon\widetilde{L}\to L$ of the links of the strata of $X$ satisfying:
\begin{enumerate}

\item the restriction $\LL\colon\LL^{-1}(X\backslash\Sigma)\longrightarrow X\backslash\Sigma$ is a smooth trivial finite covermap;

\item $\LL^{-1}(\Sigma)$ is covered by unfoldable charts $\alpha$, i.e., we have a commutative diagram
\begin{equation*}
\xymatrix{
U\times\widetilde{L}\times\R\ar[r]^(0.65){\widetilde{\alpha}}\ar[d]^{c}&  \widetilde{X}\ar[d]^{\LL}\\
U\times c(L)\ar[r]^(0.65){\alpha} &   X}
\end{equation*}
where $\widetilde{\alpha}$ is a diffeomorphism onto $\LL^{-1}(\mathop{{\rm Im}}(\alpha))$ and the left vertical map is given by $c(u,x,t)=(u,[\LL_L(x),|t|])$.
\end{enumerate}

\prg {\bf Perversities and intersection cohomology  (\cite[sec. 3]{S2})}

Recall that if $\pi\colon M\to B$ is a surjective submersion, the {\em perverse degree} $\|\omega\|_B$ of a differential form $\omega\in\Omega^*(M)$ is the smallest integer $m$ such that $i_{\xi_0}\dots i_{\xi_m}\omega =0$
for every collection of vector fields $\xi_0,\dots,\xi_m$ tangent to the fibers of $\pi$. By convention, $\|0\|_B=-\infty$.

Denote by  $\mathcal{S}^{\text{sing}}$ the set of singular strata of the stratified pseudomanifold $X$.  Recall that a {\em perversity} $\overline{p}$ in $X$ is a map $\overline{p}\colon \mathcal{S}^{\text{sing}}\to\Z$. We denote by $\mathcal{P}_X$ the set of all perversities of $X$. Fix an unfolding $\LL\colon\widetilde{X}\to X$. We say that a form $\omega\in\Omega^*(X\backslash\Sigma)$ is {\em liftable} if there exists a form $\widetilde{\omega}\in\Omega^*(\widetilde{X})$ such that $\LL^*(\omega)=\widetilde{\omega}$. The algebra of liftable forms of $X$ is denoted by $\Pi^*(X)$. Given a perversity $\overline{p}$ in $X$, recall that the cohomology of the complex
$$
\Omega_{\overline{p}}^*(X)=\{\omega\in\Pi^*(X):\ \|\omega\|_S\le\overline{p}(S)\ \text{and}\  \|d\omega\|_S\le\overline{p}(S),\ \forall S\in\mathcal{S}^{\text{sing}}\}
$$
is the $\overline{p}$-{\em intersection cohomology} of $X$, and is denoted by $\IH^*_{\overline{p}}(X)$.

\prg{\bf Modelled action  (\cite[sec. 4]{Pa}, \cite[sec. 1.1]{PS})}

Under some assumptions, the orbit space of an action of the circle on a
stratified pseudomanifold is also a stratified pseudomanifold, which is called  $\sbat$-pseudomanifold in
\cite[sec. 4]{Po}. In this work we shall use a variant of this concept: modelled actions of $\sbat$ on unfolded pseudomanifolds.
   We list  below the main properties   of a modelled action
  $\Phi \colon \sbat \times X \TO X$ of the circle  $\sbat$ on an unfolded
  pseudomanifold $X$. We denote by  $\pi \colon X \to B$ the canonical projection onto the orbit space $B =X/\sbat$.

  \Am {\em The isotropy subgroup $\sbat_x$ is the same for each $x \in S$.
  It will be denoted by $\sbat_S$}.

  \am {\em For each regular stratum $R$ we have $\sbat_{R} = \{ 1 \}$}.

  \am {\em For each  singular stratum $S$  with $\sbat_{S}= \sbat$, the
  action $\Phi$ induces a modelled action $\Phi_{L_S}
  \colon
  \sbat \times L_S \to L_S$, where $L_S$ is the link of $S$}.

  \am {\em The orbit space $B$ is an unfolded
  pseudomanifold,
  relatively to the stratification ${\cal S}_B = \{ \pi(S) \ / \ S \in {\cal
  S}_X\}$, and the projection  $\pi \colon X \to B $ is an unfolded morphism}.

  \am {\em The assignment $S \mapsto \pi (S)$ induces  the bijection $\pi_{\mathcal S}   \colon {\cal S}_X \to
  {\cal
  S}_B$.}

\noindent The action $\Phi$ may induce
  two kind of   strata of $X$:
  \begin{itemize}
 \item  a stratum $S$ is  {\em mobile} when $\tres{\S}{1}{S}$, the isotropy subgroup of any point of $S$,  is finite
  and
  \item   a stratum $S$ is is  {\em fixed} when $\tres{\S}{1}{S}$, the isotropy subgroup of any point of $S$, is $
    \sbat$.
    \end{itemize}
    Recall that the regular stratum $R$ is mobile with $\tres{\S}{1}{R} = \{1\}$. In this work, we need the refinement of
    fixed strata introduced in \cite[sec. 5.6]{Pa}
      \begin{itemize}
 \item a fixed stratum $S$ is {\em perverse\footnote{See \cite[sec. 1.2]{PS} for some examples. Notice that, as G. Friedman pointed out in \cite{F}
, there is a misprint in  \cite[sec. 1.1]{PS} in the definition of perverse stratum: it should be  $\hiru{H}{*}{L_S \menos \Sigma_{L_S}} \not=
 \hiru{H}{*}{(L_S \menos \Sigma_{L_S})/\sbat}
  \otimes
  \hiru{H}{*}{\sbat}$, where $\Sigma_{L_S}$ is the singular part of the link $L_{S}$. That definition is equivalent to the one we give above.}} when the Euler class of the action $\Phi_{\ib{L}{S}} \colon  \sbat \times \ib{L}{S} \to \ib{L}{S}$ does not vanish,  where $\ib{L}{S}$ is the link of $S$.
    \end{itemize}

\prg {\bf Gysin sequence  (\cite[sec. 6]{Pa}, \cite[sec. 1.3]{PS})\footnote{The Gysin sequence for intersection cohomology has been
constructed in \cite[sec. 6]{Pa}. In this article we use the notations of
\cite[sec. 1.3]{PS}. Notice that, as G. Friedman pointed out in \cite{F}, in  the definition of  $\lau{\mathcal{G}}{*}{\per{p}}{B}$ given in \cite[sec. 1.3]{PS}, the degree  should be shifted by 1, so that $\lau{\mathcal{G}}{*}{\per{p}}{B} \subset \lau{\Om}{*}{\per{p} - \per{x}}{B}$.}}

Fix a modelled action of $\sbat$ on $X$.
 Recall that the complex of {\em invariant $\overline{p}$-forms} computes the $\overline{p}$-intersection cohomology of $X$. This complex can be described in terms of basic data as follows: consider
 the graded
 complex
 \be
\label{inv}
\lau{\BOm}{*}{\per{p}}{X}= \left\{ (\alpha,\beta) \in
\lau{\Pi}{*}{}{B} \oplus \lau{\Om}{*-1}{\per{p}-\per{x}}{B}  \Big/
 \left[
 \begin{array}{l}
 ||\alpha||_{\pi(S)} \leq \per{p}(S) \\[,3cm] ||d\alpha + (-1)^{|\beta|}
 \beta \wedge \epsilon||_{\pi(S)} \leq \per{p}(S)
 \end{array}
 \right]
 \hbox{ if } S \in \mathcal{S}^{^{sing}}_{X} \right\} \ee endowed with
 the differential $ D(\alpha,\beta) = (d\alpha + (-1)^{|\beta|} \beta
 \wedge \epsilon , d\beta) $.  Here $| -|$ stands for the degree of the
 form, $\epsilon \in {\Pi}^{2}(B)$ is an Euler form (i.e., $\epsilon=d\chii$ for a characteristic form $\chii$ of the action) and $\per{x}$ is the characteristic perversity defined by
 $
\per{x}(\pi(S)) = \left\{
\begin{array}{ll}
1& \hbox{if $S$ is a fixed stratum }\\
0 & \hbox{if $S$ is a mobile stratum.}
\end{array}
\right.
$

\noindent The assignment $(\alpha,\beta) \mapsto \pi^{*}\alpha +
 \pi^{*}\beta \wedge \chi$ establishes a differential graded
 isomorphism between $\lau{\BOm}{*}{\per{p}}{X}$ and $ \left(
 \lau{\Omega}{*}{\per{p}} {X}\right)^{\sbat}$. This gives rise to the long exact {\em Gysin sequence}:
$$
 \cdots
 \TO
 \lau{\IH}{i+1}{\per{p}}{X} \stackrel{\ib{\oint}{\per{p}}}{\TO}
 \coho{H}{i}{\lau{{\cal G}}{*}{\per{p}}{B}}
 \stackrel{\ib{\eub}{\per{p}}}{\TO}
 \lau{\IH}{i+2}{\per{p}}{B}
 \stackrel{\ib{\pi}{\per{p}}}{\TO}
 \lau{\IH}{i+2}{\per{p}}{X}
 \TO
 \dots,
 $$
where the {\em Gysin term} $\lau{{\cal G}}{*}{\per{p}}{B}$ is the
differential
complex
 $$
\left\{
\beta \in \lau{\Om}{*-1}{\per{p} - \per{x}}{B} \ \Big/ \ \exists \alpha \in
\lau{\Pi}{*}{}{B} \hbox{ with }
\left[
\begin{array}{l}
||\alpha||_{\pi(S)} \leq \per{p}(S)\hbox { and } \\[,3cm] ||d\alpha +
(-1)^{|\beta|}\beta \wedge \epsilon||_{\pi(S)} \leq \per{p}(S)
\end{array}
\right]
\hbox{ if } S \in {\cal S}^{^{sing}}_X
\right\}.
$$
Recall that the Euler perversity $\per{e}$ is defined by
  $
    \per{e}(S) =
    \left\{
    \begin{array}{ll}
    0 & \hbox{when $S$ mobile stratum}, \\
1 & \hbox{when $S$ not perverse fixed stratum}\\
2 & \hbox{when $S$ perverse stratum.}
    \end{array}
    \right.
    $

\noindent So, the Euler class $e =  [\epsilon]$ belongs to $\lau{\IH}{2}{\per{e}}{B}$.

\prg {\bf Perverse algebras (\cite[subsection 2.2]{PS})}

A {\em perverse set} is a triple $(\mathcal{P},+,\leq)$ where
$(\mathcal{P},+)$ is an abelian semi-group with unit element
$\per{0}$ and $(\mathcal{P},\leq)$ is a partially ordered set  such that $\le$ and $+$ are compatible. Notice that the set of all perversities of an unfolded pseudomanifold $\mathcal{P}_X$ is a perverse set.

Recall that a {\em differential graded commutative (dgc, for short) perverse algebra} (or simply a perverse algebra) is a quadruple
 $\mbox{\boldmath $E$} = (E,\iota,\wedge,d)$ where
 \begin{itemize}
     \item[-] $E =\displaystyle{\bigoplus_{\per{p} \in \mathcal{P}}}\ib{E}{\per{p}}$
     where each $\ib{E}{\per{p}}$ is a graded (over $\Z$) vector space,

    \item[-]
    $\iota = \left\{ \ib{\iota}{\per{p},\per{q}} \colon \ib{E}{\per{p}} \to
    \ib{E}{\per{q}} \ / \ \per{p} \leq \per{q} \right\}$ is a family of
 graded linear morphisms,    and

  \item[-] $(E,d,\wedge)$ is a dgc algebra,
    \end{itemize}
    verifying
$$
\begin{array}{lll}
i) \ \ib{\iota}{\per{p},\per{p}} = \Ide \hspace{1cm}&

ii) \  \ib{\iota}{\per{q} ,\per{r}} \rondp \ib{\iota}{\per{p} ,\per{q}} =
      \ib{\iota}{\per{p}, \per{r}} &

iii) \ \wedge \left( \ib{E}{\per{p}} \times
	   \ib{E}{\per{p}'}\right) \subset \ib{E}{\per{p} +\per{p}'}\\[,3cm]
	
      iv) \ d\left( \ib{E}{\per{p}}\right) \subset \ib{E}{\per{p}}
&

  v) \ \ib{\iota}{\per{p}+\per{p}',\per{q}+\per{q}'}(a \wedge
      a') = \ib{\iota}{\per{p},\per{q}}(a) \wedge
      \ib{\iota}{\per{p}',\per{q}'}(a')  \hspace{1cm}
&
vi) \ d \rondp \ib{\iota}{\per{p},\per{q}} =
      \ib{\iota}{\per{p},\per{q}} \rondp d  \\
\end{array}
$$
 \nt Here, $\per{p} \leq \per{q} \leq \per{r}$, $\per{p}' \leq \per{q}'$, $a \in
 \ib{E}{\per{p}}$ and $a'\in \ib{E}{\per{p}'}.
$ \smallskip

For example, associated to a modelled action, the following dgc algebras have the structure of perverse algebras:
$
\lau{\mbox{\boldmath $\Om$}}{}{}{X} = {\displaystyle \bigoplus_{\per{p} \in \mathcal{P}_{X}}}
\lau{\Om}{}{\per{p}}{X}
$ and
$
{\displaystyle
\lau{\IH}{}{}{X} = \bigoplus_{\per{p} \in \mathcal{P}_{X}}}
\lau{\IH}{}{\per{p}}{X}
$.

\begin{remark}\label{minimal}
Notice that the intersection cohomology relatively to a fixed perversity $\per{p}$ is not an algebra, due to property $iii)$, that is,  $\bigwedge\colon\lau{\IH}{i}{\per{p}}{B}\times \lau{\IH}{j}{\per{q}}{B}\longrightarrow\lau{\IH}{i+j}{\per{p}+\per{q}}{B}$.
 With perverse algebras, we recover the algebra structure, as we consider all perversities together. The category of these objects has been studied recently in \cite{Hovey}. The perverse algebra structure has been used, too, in \cite{tanre} to extend the theory of minimal models to the context of intersection cohomology. We expect to obtain analogous results in the context of the present article.
\end{remark}

\begin{remark}
For some specific contexts such as pseudomanifolds arising from complex algebraic varieties, it would be more natural to work just with the middle perversity $\per{m}$. Nevertheless, as follows from the previous remark, multiplication by the Euler class $e\in\lau{\IH}{2}{\per{e}}{B}$ does not define an endomorphism of $\lau{\IH}{*}{\per{m}}{B}$ unless all strata are mobile, but a homomorphism $\lau{\IH}{*}{\per{m}}{B}\to \lau{\IH}{*+2}{\per{m}+\per{e}}{B}$ instead. So, to work with this homomorphism in a suitable category, we need to work with all possible perversities together, which leads us to perverse algebras.
\end{remark}

\section{Equivariant intersection cohomology}

We introduce in this section the
equivariant intersection cohomology \cite{Bry,J} of a modelled action
\cite{Pa}. For the rest of this work, we fix a modelled action
$\Phi \colon \sbat \times X \to X$. We denote by $B$ the orbit space $X/\sbat$.

\prg {\bf Equivariant intersection cohomology}. We fix $\per{p}$ a perversity of $X$.
 As $\sbat$ is connected and compact, the cohomology of the subcomplex of $\sbat$-invariant forms
$\lau{\BOm}{*}{\per{p}}{X}$ is $\lau{\IH}{*}{\per{p}}{X}$.

 Recall that the classifying space
of $\sbat$ is
just $\C\P^\infty$  whose cohomology is the free dgc algebra $\Lambda \eut$
where $|\eut| =2$ and $d\eut =0$. The {\em equivariant intersection
cohomology}
$\lau{\IH}{*}{\per{p},\sbat}{X}$,  relatively to the perversity $\per{p}$, is the cohomology of the complex
$
\left(
\lau{\BOm}{*}{\per{p}}{X} \otimes \Lambda \eut, \nabla
\right),
$
where
$\nabla$ is defined linearly from
\begin{equation*}
\nabla ((\alpha,\beta) \otimes \eut^n )=
D(\alpha,\beta) \otimes \eut^n
+
(-1)^{|\beta|}(\beta,0) \otimes \eut^{n+1}.
\end{equation*}
The equivariant intersection cohomology generalizes the usual equivariant
cohomology since  $\lau{\IH}{*}{\per{0},\sbat}{X}= \lau{H}{*}{\sbat}{X}$ when $X$ is
normal.

\prg{\bf $\Lambda \eut$-perverse algebras}. We have introduced in
\cite[sec. 2]{PS} the notion of perverse algebra, perverse morphism and
perverse isomorphism.
The coefficient ring of these objects  is $\R$. When we replace
this ring by $\Lambda \eut$, we get the notions of {\em $\Lambda \eut$-perverse
algebra}, {\em $\Lambda \eut$-perverse morphism} and
{\em $\Lambda \eut$-perverse
isomorphism}.
In this work, we deal with the following  examples:

\begin{itemize}
    \item[+] The quadruple
    $\displaystyle
 \lau{\mbox{\boldmath $\Om$}}{}{\sbat}{X} = \left(
 \bigoplus_{\per{p} \in \mathcal{P}_{X}}
  \lau{\BOm}{}{\per{p}}{X} \otimes \Lambda \eut ,\iota,\wedge,\nabla\right)
 $
 is a  $\Lambda \eut$-perverse algebra. Here, the $\Lambda \eut$-structure is given by:
$
\eut \cdot \left((\alpha,\beta) \otimes \eut^n \right)
=
 (\alpha,\beta) \otimes \eut^{n+1}.
$
 \item[+] Its cohomology
     $
     \displaystyle
 \lau{\mbox{\boldmath $\IH$}}{ }{\sbat}{X} = \left(
 \bigoplus_{\per{p} \in \mathcal{P}_{X}}
 \lau{\IH}{}{\per{p},\sbat}{X} ,\iota,\wedge,0\right)
 $
 is the the {\em equivariant intersection cohomology algebra} which
 is a $\Lambda \eut$-perverse algebra.
 \item[+]
  The operator
$\mbox{\boldmath
$\pi'$}\colon  \hiru{\mbox{\boldmath $\IH$}}{}{B} \otimes
\Lambda \eut \to
\lau{\mbox{\boldmath $\IH$}}{}{\sbat}{X}$, defined by  $\ib{\pi'}{\per{p}}\left([\alpha ]\otimes \eut^n\right)=
[(\alpha,0) \otimes \eut^n]$, is a $\Lambda \eut$-perverse morphism.

\end{itemize}

\begin{remark}
The  $\Lambda \eut$-perverse morphism $\pi'$ suggests that the equivariant perverse minimal model of $X$ (see Remark \ref{minimal}) may be computed by the mimimal model of $B$, as it happens in the case without perversities.
\end{remark}

\prg {\bf Equivariant  Gysin sequence}. The main tool we use  for the
classification of modelled actions is the Gysin sequence  we construct now. Fix $\per{p}$  a perversity of $X$. Consider the short
exact sequence
\begin{equation*}
0 \to
\left(
\lau{\Om}{*}{\per{p}}{B} \otimes \Lambda \eut,d \otimes 1
\right)
\stackrel{\ib{\pi'}{\per{p}}}{\TO}
\left(
\lau{\BOm}{*}{\per{p}}{X} \otimes \Lambda \eut,\nabla
\right)
\stackrel{\ib{\oint'}{\per{p}}}{\TO}
\left(
\lau{\mathcal{G}}{*-1}{\per{p}}{B} \otimes \Lambda \eut,d \otimes 1
\right)
\to
0,
\end{equation*}
where
$\ib{\pi'}{\per{p}}(\alpha \otimes \eut^n) = (\alpha,0) \otimes \eut^n$
and
$\ib{\oint'}{\per{p}}(\alpha,\beta) \otimes \eut^n = \beta \otimes \eut^n$.
Each term is a
differential complex and a $\Lambda \eut$-module with the natural
structure. Moreover, the maps $\ib{\pi'}{\per{p}}$ and
$\ib{\oint'}{\per{p}}$  preserve these  structures.
The {\em equivariant  Gysin sequence} is the induced long exact sequence
$$
\cdots \to
\left[\lau{\IH}{*}{\per{p}}{B} \otimes \Lambda \eut \right]^i
\stackrel{\ib{\pi'}{\per{p}}}{\TO}
\lau{\IH}{i}{\per{p},\sbat}{X}
\stackrel{\ib{\oint'}{\per{p}}}{\TO}
\left[\hiru{H}{*}{\lau{\mathcal{G}}{\cdot}{\per{p}}{B}} \otimes \Lambda
\eut \right]^{i-1}
\stackrel{\ib{\delta}{\per{p}}}{\TO}
\left[\lau{\IH}{*}{\per{p}}{B} \otimes \Lambda \eut \right]^{i+1}
\to \cdots
$$
Here,
$
\ib{\delta}{\per{p}}
\left([\beta ]\otimes \eut^n\right)
=
[d\alpha + (-1)^{|\beta|} \beta \wedge \epsilon ]\otimes \eut^n
+
(-1)^{|\beta|}\tres{\iota}{B}{\per{p}-
\per{x}, \per{p}}[\beta]\otimes \eut^{n+1}.
$
For short, we shall write $(-1)^{|\beta|}\tres{\iota}{B}{\per{p}-
\per{x}, \per{p}}{} = \ib{I}{\per{p}}$. The connecting morphism
becomes $\ib{\delta}{\per{p}}=
\ib{\eub}{\per{p}} \otimes 1+ \ib{I}{\per{p}} \otimes \eut  .
$
Notice that the equivariant Gysin sequence permits us to obtain the
equivariant intersection cohomology  of $X$ in terms of
basic\footnote{Of the orbit space $B$.} data.

\section{Classification of modelled actions}

In this section, we prove that $B$ and the Euler class determine the equivariant intersection cohomology of
$X$.

\prg {\bf Fixing the orbit space}\footnote{See (cf. \cite[sec. 3.1]{PS}) for details.}.
Consider $\Phi_{1} \colon \sbat \times X_{1}\to X_{1}$ and $\Phi_{2}
\colon \sbat \times X_{2}\to X_{2}$ two modelled actions and write
$B_{1}$ and $B_{2}$ the corresponding orbit spaces.

An  unfolded isomorphism $f \colon B_{1} \to B_{2}$  is {\em optimal} when it preserves the nature
of the strata. In this case, the two Euler perversities are
equal: $\per{e}_{1}(\pi_{1}(S)) = \per{e}_{2}(f(\pi_{1}(S)))$ for each
singular stratum $S$ of $X_1$.  We shall write
$\per{e}$ for this {\em Euler perversity}.
Now we can compare the two Euler classes
$e_{1}\in \lau{\IH}{2}{\per{e}}{B_{1}}$ and $e_{2} \in \lau{\IH}{2}{\per{e}}{B_{2}}$.
We shall say that $e_{1}$ and $e_{2}$ are $f$-{\em related} if
$\tres{f}{*}{\per{e}}e_{2} =  e_{1}$.

\bp
\label{aritz}
Let $X_{1}$, $X_{2}$ be two connected normal unfolded pseudomanifolds.
Consider two modelled actions $\Phi_1 \colon \sbat \times X_1 \to X_1$
and $\Phi_2 \colon \sbat \times X _2\to X_2$.  Let us suppose that
there exists an unfolded isomorphism $f \colon B_{1} \to B_{2}$
between the associated orbit spaces.  Then, the two following
statements are equivalent:

\Zati The isomorphism $f$ is optimal and the Euler classes $e_{1}$
and $e_{2}$ are $f$-related.

\zati There exists a $\Lambda \eut$-perverse isomorphism
$\mbox{\boldmath $G$} \colon \lau{\mbox{\boldmath $\IH$}}{}{\sbat}{X_{2}}\to
\lau{\mbox{\boldmath $\IH$}}{}{\sbat}{X_{1}}$
verifying
$
\mbox{\boldmath $G$} \rondp \mbox{\boldmath $\pi'_{2}$} = \mbox{\boldmath
$\pi'_{1}$} \rondp (\mbox{\boldmath $f$} \otimes1).$
\ep
\pro We proceed in two steps.

\medskip

\fbox{$(a) \Rightarrow (b)$}   Since $[\tres{f}{*}{\per{e}}\epsilon_2] =\tres{f}{*}{\per{e}}e_{2}   =e_{1} =
[\epsilon_1] $ then there exists $\gamma \in
\lau{\Om}{1}{\per{e}}{B_{2}}$ with $\tres{f}{*}{\per{e}}\epsilon_2 =
\epsilon_1 - d(\tres{f}{*}{\per{e}}\gamma)$.
Using this map $\gamma$, we construct for each perversity $\per{p}$ the map
$\ib{G}{\per{p}} =  \ib{F}{\per{p}} \otimes 1 \colon \lau{\Omega}{*}
{\per{p},\sbat} {X_2}\TO
\lau{\Omega}{*}{\per{p},\sbat} {X_1}$ (cf. \cite[Proposition 3.2]{PS}).
Since $\mbox{\boldmath $F$} = \{ \ib{F}{\per{p}} \}$ is a perverse
isomorphism, then  $\mbox{\boldmath $G$} = \{ \ib{G}{\per{p}} \} \colon
\lau{\mbox{\boldmath $\Om$}}{}{\sbat}{X_{2}}\to \lau{\mbox{\boldmath
$\Om$}}{}{\sbat}{X_{1}}$ is a $\Lambda \eut$-perverse isomorphism.
The equality $
\mbox{\boldmath $G$} \rondp \mbox{\boldmath $\pi'_{2}$} = \mbox{\boldmath
$\pi'_{1}$} \rondp (\mbox{\boldmath $f$} \times 1 )$
comes from:
$$
\ib{G}{\per{p}} \left( \ib{\pi'}{2,\per{p}}
\left(\alpha \otimes \eut^n\right)\right)
=
\ib{F}{\per{p}} (\alpha, 0) \otimes  \eut^n
=
(\ib{f}{\per{p}}(\alpha ) ,0) \otimes \eut^n=
\ib{\pi'}{1,\per{p}} \left( \ib{f}{\per{p}}(\alpha ) \otimes \eut^n\right),
$$
where  $\alpha \in
\lau{\Om}{*}{\per{p}}{B_{2}}$.

\medskip

\fbox{$(b) \Rightarrow (a)$}
Consider now the  equivariant Gysin sequences associated to the actions $\Phi_1$
and
$\Phi_2$.
The two Gysin terms are written {$_{_1}\mathcal{G}$} and
{$_{_2}\mathcal{G}$} respectively.
Since
$
\ib{G}{\per{e}_{2}} \rondp \ib{\pi'}{2,\per{e}_{2}}
=
\ib{\pi'}{1,\per{e}_{2}} \rondp (\ib{f}{\per{e}_{2}} \times 1)
$
we can construct the commutative diagram
$$
\begin{CD}
\lau{\IH}{1}{\per{e}_2}{B_{2}}
@>\ib{\pi'}{2,\per{e}_2}>>\lau{\IH}{1}{\per{e}_2, \sbat}{X_{2}}
@>\ib{\oint'}{2,\per{e}_2}>>
\hiru{H}{0}{\lau{_{_1}\mathcal{G}}{*}{\per{e}_2}{B_{2}}  }
@>\ib{\delta'}{2,\per{e}_2}>>
\left[\lau{\IH}{*}{\per{e}_2}{B_{2}} \otimes \Lambda \eut\right]^2
@>\ib{\pi'}{2,\per{e}_2}>>\lau{\IH}{2}{\per{e}_2,\sbat}{X_{2}} \\
    @VV\ib{f}{\per{e}_{2}}V @V\ib{G}{\per{e}_2} VV @V\ell VV
@V V\ib{f}{\per{e}_{2}}\otimes 1 V @V\ib{G}{\per{e}_2}  VV  \\
\lau{\IH}{1}{\per{e}_2}{B_{1}}
@>\ib{\pi'}{1,\per{e}_2}>>\lau{\IH}{1}{\per{e}_2,\sbat}{X_{1}}
@>\ib{\oint'}{1,\per{e}_2}>>
\hiru{H}{0}{\lau{_{_2}\mathcal{G}}{*}{\per{e}_2}{B_{1}} }
@>\ib{\delta'}{1,\per{e}_2}>>
\left[\lau{\IH}{*}{\per{e}_2}{B_{1}} \otimes \Lambda \eut\right]^2
@>\ib{\pi'}{1,\per{e}_2}>>\lau{\IH}{2}{\per{e}_2,\sbat}{X_{1}} ,\\
    \end{CD}
$$
where $\ell  \colon
\hiru{H}{0}{\lau{_{_1}\mathcal{G}}{*}{\per{e}_1}{B}} \to
\hiru{H}{0}{\lau{_{_2}\mathcal{G}}{*}{\per{e}_1}{B}}$
 is an isomorphism. Following the proof of  \cite[Proposition 3.2]{PS}
 we conclude that  the isomorphism $f$ is optimal, the operator
 $\ell$ is the multiplication by  a number
 $\lambda \in \R \menos \{
0\}$ and
$\ib{f}{\per{e}}  e_2 = \lambda \cdot e_1$.
 Finally, the commutativity  $\left( \ib{f}{\per{e}}\otimes 1
 \right)
 \rondp \ib{\delta'}{2,\per{e}_2}
 = \ib{\delta'}{1,\per{e}_2} \rondp \ell$ gives that $\lambda =1$ and
 therefore  the Euler classes $e_{1}$ and $e_{2}$ are $f$-related.\qed

\section{The basic spectral sequence}

The Leray spectral sequence considered by
Borel for the usual equivariant cohomology has been extended to the perverse framework
in \cite{Bry}. It converges to $\lau{\IH}{*}{\per{p},\sbat}{X}$ and its
second term is $\lau{\IH}{*}{\per{p}}{X} \otimes  \Lambda \eut$.

 We construct another
spectral
sequence converging to $\lau{\IH}{*}{\per{p},\sbat}{X}$ whose
third term  is described in terms of
 $B$.
	It is the {\em basic spectral sequence}.
First of all, we present an auxiliary complex.

\prg {\bf The co-Gysin complex}. The third term of the spectral sequence
is described in terms of  the  {\em co-Gysin complex}\footnote{An element of $\lau{\mathcal{K}}{*}{\per{p}}{B}$ is written $\overline\alpha$ where $\alpha \in \lau{\Om}{*}{\per{p}}{B}.$}
$
\displaystyle
\lau{\mathcal{K}}{*}{\per{p}}{B} =
\frac{\lau{\Om}{*}{\per{p}}{B}}{\lau{\mathcal{G}}{*}{\per{p}}{B}}$.
  It fits into the long exact sequence
$$
\cdots
 \to
\hiru{H}{i}{\lau{\mathcal{G}}{*}{\per{p}}{B}}
\stackrel{\ib{I}{\per{p}}}{\TO}
\lau{\IH}{i}{\per{p}}{B}
\stackrel{\ib{P}{\per{p}}}{\TO}
\hiru{H}{i}{\lau{\mathcal{K}}{*}{\per{p}}{B}}
\stackrel{\ib{\partial}{\per{p}}}{\TO}
\hiru{H}{i+1}{\lau{\mathcal{G}}{*}{\per{p}}{B}}
\to
\cdots.
$$
Here, $\ib{I}{\per{p}}[\alpha] = [\alpha]$,
$\ib{\partial}{\per{p}}[\overline\alpha] =
[d\alpha]$ and $\ib{P}{\per{p}}[\alpha] = [\overline\alpha]$.
Now, we can describe the basic spectral sequence.

\bp
\label{spec}
Consider a modelled action $\Phi \colon \sbat \times X \to X$ and fix a
perversity $\per{p}$. There exists a first quadrant spectral sequence
$\left\{ \left( \ib{E}{\per{p},r} , \ib{d}{\per{p},r} \right)\right\}_{r\geq 0}$ converging to the
equivariant intersection cohomology $\lau{\IH}{*}{\per{p},\sbat}{X}$ such that

\Zati  $\tres{E}{i,j}{\per{p},r}  = 0$ if $j$ is an odd number and $r
\geq 1$;

\zati $\tres{E}{i,2j}{\per{p},2s} =\tres{E}{i,2j}{\per{p},2s+1} $ if $s\geq 1$;

\zati the second and third terms are
$
\tres{E}{i,2j}{\per{p},3} =
\tres{E}{i,2j}{\per{p},2} =
\left\{
\begin{array}{cl}
\lau{\IH}{i}{\per{p}}{B}
& \hbox{if } j=0 \\[,2cm]
\hiru{H}{i}{\lau{\mathcal{K}}{*}{\per{p}}{B}} \otimes
\R \cdot\eut^j & \hbox{if } j>0;
\end{array}
\right.
$

\zati the third differential
$
\ib{d}{\per{p},3} \colon
\tres{E}{i,2j}{\per{p},3} \TO
\tres{E}{i+3,2j-2}{\per{p},3}
$
is
\ $
\ib{d}{\per{p},3} (w \otimes \eut^j) =
\left\{
\begin{array}{cl}
\left(\ib{\eub}{\per{p}}\rondp \ib{\partial}{\per{p}}\right)(w)
&\hbox{if } j=1
\\[,2cm]
\left(\ib{P}{\per{p}}\rondp\ib{\eub}{\per{p}}\rondp\ib{\partial}{\per{p}}\right)(w) \otimes
\eut^{j-1}
&\hbox{if } j\geq 2.
\end{array}
\right.
$
\ep
\pro
Consider the filtration
$
\cdots \subset
\bi{F}{i}{\lau{\Om}{*}{\per{p},\sbat}{X}}
\subset
\bi{F}{i-1}{\lau{\Om}{*}{\per{p},\sbat}{X}}
\subset
\cdots
\subset
\bi{F}{0}{\lau{\Om}{*}{\per{p},\sbat}{X}}
=
\lau{\Om}{*}{\per{p},\sbat}{X}
$
defined by
$
\bi{F}{i}{\lau{\Om}{*}{\per{p},\sbat}{X}} =
\{ \om \in \lau{\BOm}{\geq i}{\per{p}}{X} \otimes \Lambda \eut \ / \
d\om \in \lau{\BOm}{\geq i}{\per{p}}{X} \otimes \Lambda \eut \}.
$
That is,
\begin{eqnarray*}
\bi{F}{i}{\lau{\Om}{i+2j}{\per{p},\sbat}{X}} & =&
 \left(\lau{\Om}{i}{\per{p}}{B} \oplus \{ 0 \} \right)\otimes \R \cdot \eut^j
\oplus \bigoplus_{k=1}^{j}
\lau{\BOm}{i+2k}{\per{p}}{X} \otimes \R \cdot \eut^{j-k}, \hbox{ and}
\\
\bi{F}{i}{\lau{\Om}{i+2j+1}{\per{p},\sbat}{X}}&  =&	
 \lau{\BOm}{i+1}{\per{p}}{X} \otimes \R \cdot \eut^j
\oplus \bigoplus_{k=1}^{j}
\lau{\BOm}{i+1+2k}{\per{p}}{X} \otimes \R \cdot \eut^{j-k}
\end{eqnarray*}
which verifies
$
\nabla \left( \bi{F}{i}{\lau{\Om}{*}{\per{p},\sbat}{X}} \right) \subset
\bi{F}{i}{\lau{\Om}{*+1}{\per{p},\sbat}{X}} .
$
Following the standard procedure (see for example \cite{MC})
one constructs a spectral sequence $\left\{ \left( \ib{E}{\per{p},r} ,
\ib{d}{\per{p},r} \right)\right\}$ converging to the
equivariant cohomology $\lau{\IH}{*}{\per{p},\sbat}{X}$. We have

$$
\tres{E}{i,2j}{\per{p},0}
=
\left(\lau{\Om}{i}{\per{p}}{B} \oplus \{ 0 \}\right)\otimes \eut^{j}
\ \ \ \hbox{ and } \ \ \
\tres{E}{i,2j+1}{\per{p},0}
=
\frac{\lau{\BOm}{i+1}{\per{p}}{X}}{\lau{\Om}{i+1}{\per{p}}{B} \oplus \{ 0
\}} \otimes  \eut^{j}.
$$
The differential $\ib{d}{\per{p},0} \colon \tres{E}{i,2j}{\per{p},0}  \to
\tres{E}{i,2j+1}{\per{p},0} $
is zero, and the differential
$\ib{d}{\per{p},0} \colon \tres{E}{i,2j+1}{\per{p},0}  \to
\tres{E}{i,2j+2}{\per{p},0} $
is given by
$\ib{d}{\per{p},r} \left(\overline{(\alpha,\beta)}\otimes \eut^j \right)= (\beta,0)
\otimes \eut^{j+1}.
$
We conclude that
$
\tres{E}{i,j'}{\per{p},1}
=
\left\{
\begin{array}{cl}
\lau{\Om}{i}{\per{p}}{B} & \hbox{if $j'=0$ }\\[,2cm]
\lau{\mathcal{K}}{i}{\per{p}}{B} \otimes \R \cdot\eut^{j} & \hbox{if $j'= 2 j>0$}\\[,2cm]
0&\hbox{if $j'$ is
odd}
\end{array}
\right.
$.
This gives (a),  $\ib{d}{\per{p},2s}= 0$ if $s\geq 1$ and (b).

The first differential
$\ib{d}{\per{p},1} \colon \tres{E}{i,2j}{\per{p},1}  \to
\tres{E}{i+1,2j}{\per{p},1} $
is given by
$
\ib{d}{\per{p},1} (\alpha) = d \alpha$ and $
\ib{d}{\per{p},1}\left(\overline{\alpha}\otimes \eut^j \right)= d
\overline{\alpha}
\otimes \eut^j.
$
We conclude that
$
\tres{E}{i,2j}{\per{p},3} = \tres{E}{i,2j}{\per{p},2} =
\left\{
\begin{array}{cl}
\lau{\IH}{i}{\per{p}}{B}
& \hbox{if } j=0 \\[,2cm]
\hiru{H}{i}{\lau{\mathcal{K}}{*}{\per{p}}{B}} \otimes
\R \cdot\eut^j & \hbox{if } j>0.
\end{array}
\right.
$.
This gives (c).

Consider, for the computation of the third differential,
$[\overline{\alpha}] \in\hiru{H}{i}{\lau{\mathcal{K}}{*}{\per{p}}{B}}$.
So, we have that $[d\alpha] \in \hiru{H}{i+1}{\lau{\mathcal{G}}{*}{\per{p}}{B}}$
and  $\ib{\eub}{\per{p}}([d\alpha] )\in \lau{\IH}{i+3}{\per{p}}{B}$.
This gives $\ib{d}{\per{p},3} ([\overline{\alpha}] \otimes \eut)  =
\ib{\eub}{\per{p}}([d\alpha] )
=
\ib{\eub}{\per{p}}\ib{\partial}{\per{p}}([\overline{\alpha}])$.
For the general case $j\geq 2$ we have
$\ib{d}{\per{p},3} ([\overline{\alpha}] \otimes \eut^j)
=
\ib{P}{\per{p}}\ib{\eub}{\per{p}}([d\alpha] ) \otimes \eut^{j-1}
=
\ib{P}{\per{p}}\ib{\eub}{\per{p}}\ib{\partial}{\per{p}}([\overline{\alpha}] )
\otimes \eut^{j-1}.
$
This gives (d).
\qed

\prg {\bf The basic spectral sequence in the classic framework.} We
consider here the usual cohomology, that is, the case $\per{p} =
\per{0}$.
For the sake of simplicity we also suppose that $X$ is normal.

In
this context, the basic spectral sequence is a  spectral sequence
converging to $\lau{H}{*}{\sbat}{X}$ whose third term is described
in terms of the cohomology of $B$ and $F$,  the union of fixed strata.
In fact, as C. Allday pointed out to us, this spectral sequence degenerates into the
Skjelbred exact sequence (cf. \cite{Sk}).

First of all, we fix some facts.
The cohomology of the complex $\lau{\mathcal{G}}{*}{-\per{x}}{B} =
\lau{\Om}{*}{-\per{e}}{B} =
\lau{\Om}{*}{-\per{x}}{B} =
\lau{\mathcal{G}}{*}{\per{0}}{B}
$
is $\hiru{H}{*}{B,F}$ (cf. \cite{S2}).
We shall write
$\eub\colon \hiru{H}{*}{B,F} \to \hiru{H}{*+2}{B,F}$
the
map induced from $\ib{\eub}{-\per{x}}$.
The long exact sequence associated to the pair $(B,F)$ is
$
\cdots  \to \hiru{H}{i}{B,F} \stackrel{\iota}{\TO}
\hiru{H}{i}{B}
\stackrel{P}{\TO}
\hiru{H}{i}{F}
\stackrel{\partial}{\TO}
\hiru{H}{i+1}{B,F}
\to \cdots .
$

\bl
\label{bat}
Consider a modelled action $\Phi \colon \sbat \times X \to X$ where $X$ is normal. There exists a first quadrant spectral sequence
$\left\{ \left( \ib{E}{r} , \ib{d}{r} \right)\right\}_{r\geq 0}$ converging to  $\lau{H}{*}{\sbat}{X}$ such that

\Zati  $\tres{E}{i,j}{r}  = 0$ if $j$ is an odd number and $r
\geq 1$;

\zati $\tres{E}{i,2j}{2s} =\tres{E}{i,2j}{2s+1} $ if $s\geq 1$;

\zati
$
\tres{E}{i,2j}{3} =
\tres{E}{i,2j}{2} =
\left\{
\begin{array}{cl}
\hiru{H}{i}{B}
& \hbox{if } j=0 \\[,2cm]
\hiru{H}{i}{F} \otimes
\R \cdot\eut^j & \hbox{if } j>0;
\end{array}
\right.
$

\zati each $\tres{E}{*,0}{r} $ is a quotient of $\hiru{H}{*}{B}$
when $r\geq 3$;

\zati for each $s\geq 1$ and $j \ne s$, the differential
$
\ib{d}{2s+1} \colon \tres{E}{i,2j}{2s+1} \to \tres{E}{i+2s+1,2j-2s}{2s+1}
$
is 0;

\zati for each $s\geq 1$, the differential
$
\ib{d}{2s+1} \colon \tres{E}{i,2s}{2s+1}  = \hiru{H}{i}{F} \otimes \R \cdot
\eut^s\to \tres{E}{i+2s+1,0}{2s+1}
$
is induced by  $ (-1)^s \, \iota \rondp
\eub^s  \rondp \,
\partial$.
\el
\pro
Consider $\left\{\left(\ib{E}{r},\ib{d}{r}\right)\right\}_{r\geq 0}$ the spectral sequence given
by the above Proposition for $\per{p}=\per{0}$. It converges to
$\lau{\IH}{*}{\per{0},\sbat}{X}$ which is $ \lau{H}{*}{\sbat}{X}$ since $X$
is normal. Let us verify the properties.

\Zati and \zati Clear.

\zati Since $X$ is normal, then $B$ is normal and therefore
$\lau{\IH}{*}{\per{0}}{B} = \hiru{H}{*}{B}$.
The long exact sequence associated to the short exact sequence
$
0 \to \lau{\Om}{*}{-\per{x}}{B} \to
\lau{\Om}{*}{\per{0}}{B}
\to
\lau{\mathcal{K}}{*}{\per{0}}{B}
\to
0
$
becomes
$$
\cdots  \to \hiru{H}{i}{B,F} \stackrel{\iota}{\TO}
\hiru{H}{i}{B}
\stackrel{\ib{P}{\per{0}}}{\TO}
\hiru{H}{i}{\lau{\mathcal{K}}{\cdot}{\per{0}}{B}}
\stackrel{\ib{\partial}{\per{0}}}{\TO}
\hiru{H}{i+1}{B,F}
\to \cdots.
$$
So, there exists an isomorphism
$ \xi \colon \hiru{H}{*}{\lau{\mathcal{K}}{\cdot}{\per{0}}{B}}
\to
\hiru{H}{*}{F}$ with $\partial \rondp \xi= \ib{\partial}{\per{0}}$ and $\xi
\rondp \ib{P}{\per{0}} = P$. The result comes now directly from
Proposition \ref{spec}.

\zati For $r = 2s+1 \geq 3$ we have
$
\tres{E}{i,0}{2s+1}  =
\fracc{\tres{Z}{i,0}{2s+1}}{\tres{B}{i,0}{2s}}
=
\fracc{\hiru{\Om}{i}{B} \cap d^{-1}(0) }{\tres{B}{i,0}{2s}}
=
\fracc{\hiru{H}{i}{B}}{\tres{B}{i,0}{2s}\big/ d\hiru{\Om}{i-1}{B}} .
$

To prove \zati  and \zati, we proceed by induction on $s$. Taking $s=1$, we have from
 Proposition \ref{spec} and the above identifications:
$
\ib{d}{3} (w \otimes \eut^j) =
\left\{
\begin{array}{cl}
\left(\iota \rondp \eub  \rondp \partial \right)(w)
&\hbox{if } j=1
\\[,2cm]
\left(\xi \rondp \ib{P}{\per{0}}\rondp \iota \rondp \eub \rondp \partial\right)(w) \otimes
\eut^{j-1}
&\hbox{if } j\geq 2
\end{array}
\right.
=
\left\{
\begin{array}{cl}
\left(\iota \rondp \eub  \rondp \partial \right)(w)
&\hbox{if } j=1
\\[,2cm]
0
&\hbox{if } j\geq 2.
\end{array}
\right.
$.
Let us now suppose that the result is true for $s' < s$. The
case $j < s$ is straightforward by dimension reasons. Consider now $j \geq s$.
The induction hypothesis
and (b)
give the isomorphism chain
$$
\ib{\nabla}{j,s} \colon \tres{E}{i,2j}{2s+1}  =
\frac{\tres{Z}{i,2j}{2s+1}}{\tres{Z}{i+1,2j-1}{2s} + \tres{B}{i,2j}{2s}} \to
\tres{E}{i,2j}{2s} \to
\tres{E}{i,2j}{2s-1} \to
\cdots
\to
\tres{E}{i,2j}{2} \to
\hiru{H}{i}{F} \otimes \R \cdot \eut^j,
$$
with
$
\ib{\nabla}{j,s}\left(
\mbox{\boldmath $\om_j$} =
\overline{(\alpha_0,0) \otimes  \eut^j + \suma{k=1}{j}
(\alpha_k,\beta_k) \otimes
\eut^{j-k} }\right)
 =  \xi\left[\overline{\alpha_0}\right] \otimes \eut^j.
$
On the other hand, we have that the definition of the differential
$
\ib{d}{2s+1}\colon \tres{E}{i,2j}{2s+1}  =
\fracc{\tres{Z}{i,2j}{2s+1}}{\tres{Z}{i+1,2j-1}{2s} + \tres{B}{i,2j}{2s}}\TO
\fracc{\tres{Z}{i+2s+1,2j-2s}{2s+1}}{\tres{Z}{i+2s+2,2j-2s-1}{2s}
+\tres{B}{i+2s+1,2j-2s}{2s}}
$
is
$$
\ib{d}{2s+1}(\mbox{\boldmath $\om_j$} )=
\overline{(d\alpha_{s}+(-1)^{|\beta_{s} |}\beta_{s}  \wedge \epsilon+
(-1)^{|\beta_{s+1}|}\beta_{s+1} ,0) \otimes  \eut^{j-s} + \sum_{k=s+1}^{j}
(\alpha'_k,\beta'_k) \otimes
\eut^{j-k}},
$$
and therefore
$
(\ib{\nabla}{j-s,s} 	\rondp \ib{d}{2s+1})(\mbox{\boldmath $\om_j$}  )
= \xi\left[\overline{d\alpha_{s}+(-1)^{|\beta_{s} |}\beta_{s}  \wedge \epsilon +
(-1)^{|\beta_{s+1}|}\beta_{s+1}}\right]  \otimes \eut^{j-s}= 0
$
since $\beta_{s}  \wedge \epsilon \in \lau{\Om}{*}{-\per{x}}{B}$.
This implies
$
\ib{d}{2s+1}(\mbox{\boldmath $\om_j$} ) =0
$ for $j> s$.
It remains the case $j=s$. For $\mbox{\boldmath $\om_s$}  =
\ib{\nabla}{s,s}^{^{-1}}([\omega] \otimes \eut^s) \in \tres{E}{i,2s}{2s+1}  =
\fracc{\tres{Z}{i,2s}{2s+1}}{\tres{B}{i,2s}{2s}}$ we have from (d)
\begin{eqnarray*}
\ib{d}{2s+1}(\mbox{\boldmath $\om_s$} )&=&
\overline{\left[ d\alpha_s + (-1)^{|\beta_s|} \beta_s \wedge \epsilon\right]}
=
\overline{
\left[
\suma{k=0}{s} (-1)^{s-k}d\alpha_k \wedge \epsilon^{s-k}
\right]
}
=
(-1)^s \overline{\left[d\alpha_0 \wedge
\epsilon^{s}\right]
}
=\\[,2cm]
&=&
(-1)^s \overline{(\iota \rondp \eub^s  \rondp \, \ib{\partial}{\per{0}})
[\overline{\alpha_0}] }
=
(-1)^s \overline{(\iota \rondp \eub^s  \rondp \, \partial)
\xi[\overline{\alpha_0}] }
=
(-1)^s \overline{(\iota \rondp \eub^s  \rondp \, \partial)[\omega]},
\end {eqnarray*}
since
$
\suma{k=1}{s} (-1)^{s-k}d\alpha_k \wedge
\epsilon^{s-k} =
d\left(\suma{k=1}{s} \left( \alpha_k,(-1)^{|\alpha_k|}
\suma{j=1}{k-1} (-1)^{k-1-j}d\alpha_j \wedge
\epsilon^{k-1-j}\right) \otimes \eut^{s-k}\right) $ belongs to $\tres{B}
{i,0
}
{
    2s
}
.$
\qed

The particular geometry of this spectral sequence gives rise to
a {\em Gysin sequence} (in the sense of \cite{MC}).

\bp
\label{skj}
Consider a modelled action $\Phi \colon \sbat \times X \to X$ where $X$ is normal. We have the Skjelbred exact sequence
$$
\cdots
\TO
\left[ \hiru{H}{*}{F} \otimes \Lambda^{>0} \eut \right]^i
\stackrel{\mbox{\boldmath $\beta$}}{\TO}
\hiru{H}{i+1}{B}
\stackrel{\mbox{\boldmath $\alpha$}}{\TO}
\lau{\IH}{i+1}{\sbat}{X}
\stackrel{\mbox{\boldmath $\delta$}}{\TO}
\left[ \hiru{H}{*}{F} \otimes \Lambda^{>0} \eut \right]^{i+1}
\stackrel{\mbox{\boldmath $\beta$}}{\TO}
\hiru{H}{i+2}{B}
\cdots,
$$
where
\begin{itemize}
\item[-] \mbox{\boldmath $\alpha$}$[\alpha] = [(\alpha,0) \otimes 1]$;

\item[-] \mbox{\boldmath $\beta$}$
([\omega] \otimes \eut^s)=
 (-1)^s\, (\iota \rondp \eub^s  \rondp \, \partial )[\omega]$;

\item[-] \mbox{\boldmath $\delta$}$
\left[ \suma{k=0}{j} (\alpha_k,\beta_k) \otimes \eut^k\right]
=
\suma{k=1}{j} \xi[\overline{\alpha_k}] \otimes \eut^k.$
\end{itemize}
\ep
\pro
Consider the exact sequence
$\displaystyle
0
\TO
\bigoplus_{s\geq1}\tres{E}{i-2s,2s}{\infty}
\stackrel{\nabla}{\TO}
\bigoplus_{s\geq1}  \hiru{H}{i-2s}{F} \otimes \eut^s
\stackrel{\mbox{\boldmath $\beta$}}{\TO}
\hiru{H}{i+1}{B}
\stackrel{\proj}{\TO}
\tres{E}{i+1,0}{\infty}
\TO
0,
$
where   $\nabla$ is induced by $\ib{\nabla}{s,s+1}
\colon\tres{E}{i-2s,2s}{\infty} = \tres{E}{i-2s,2s}{2s+3} \to
\hiru{H}{i-2s}{F} \otimes \R \cdot \eut^s$,
and proceed as in \cite[pag. 8]{MC}.
\qed

\section{Localization}
The localization of the equivariant intersection cohomology is a cohomological theory
introduced in \cite{Bry,J}. In fact, it is a residual cohomology
since it depends on a neighborhood of the fixed point set $F$. The usual\footnote{When $X$ is a manifold, the family of strata $\mathcal{S}_X$ is reduced to the regular stratum.}
LocalizationTheorem establishes that the localization $\lau{L}{*}{\sbat}{X}$ of $\lau{H}{*}{\sbat}{X}$ is in fact $\lau{H}{*}{}{F} \otimes \R(\eut)$.
This doesn't hold for the generic case since the links of strata are
no longer spheres.

\prg {\bf Definition and properties.} Denote by $\R(\eut)$ the field of
fractions of $\Lambda
\eut$. The localization of the equivariant intersection
cohomology is
$
\lau{\IL}{\star}{\per{p},\sbat}{X} =
\lau{\IH}{*}{\per{p},\sbat}{X} \ib{\otimes}{\Lambda\euts} \R(\eut).
$
It is not a graded $\R(\eut)$-vector space  over $\Z$ but over $\Z_2$ by:
$
\lau{\IL}{\star}{\per{p},\sbat}{X} =
\lau{\IL}{even}{\per{p},\sbat}{X}
\oplus
\lau{\IL}{odd}{\per{p},\sbat}{X}
=
\lau{\IH}{even}{\per{p},\sbat}{X} \ib{\otimes}{\Lambda\euts} \R(\eut)
\oplus
\lau{\IH}{odd}{\per{p},\sbat}{X} \ib{\otimes}{\Lambda\euts} \R(\eut).
$
It verifies the following properties.

\Zati The localization $\lau{\IL}{\star}{\per{0},\sbat}{X}$ is
the usual localization $
\lau{L}{*}{\sbat}{X}$ when $X$ is  normal.

\zati A perverse algebra is a {\em perverse superalgebra} when
 the coefficient ring is $\R(\eut)$ (instead of $\R$) and it is graded over
 $\Z_{2}$ (instead of over $\Z$).
In the same manner, we introduce the notions of {\em perverse
superalgebra morphism}
and {\em perverse superalgebra isomorphism}.
The quadruple
$
\displaystyle
 \lau{\mbox{\boldmath $\IL$}}{}{\sbat}{X} = \left(
 \bigoplus_{\per{p} \in \mathcal{P}_{X}}
 \lau{\IL}{}{\per{p},\sbat}{X} ,\iota,\wedge,0\right)
 $
 is a perverse superalgebra.
On the other hand, the operator
$\mbox{\boldmath
$\pi''$}\colon  \hiru{\mbox{\boldmath $\IH$}}{}{B} \otimes
\R(\eut ) \to
\lau{\mbox{\boldmath $\IL$}}{}{\sbat}{X}$, defined by
$\ib{\pi''}{\per{p}} (b\otimes P) =
\ib{\pi'}{\per{p}}(b \otimes 1) \ib{\otimes}{\Lambda \euts} P$
, is a perverse superalgebra morphism.

\zati The localization of the equivariant intersection cohomology is a
residual cohomology: the inclusion induces a perverse superalgebra isomorphism
$\lau{\IL} {*}{\per{p},\sbat}{X} \cong \lau{\IL}
{*}{\per{p},\sbat}{U}$ where $U \subset X$ is any neighborhood of the fixed point
set $F$. In fact, $\lau{\IL} {*}{\per{p},\sbat}{X}$ can be seen as the
global sections of a sheaf $\ib{\mathcal{H}}{\per{p}}$ defined on $F$. This sheaf is {\em
constructible}, that is, locally constant on each fixed stratum $S$. Its
stalk is given by \refp{amaia}.

\zati Let us suppose that $X$ is compact. Given two complementary perversities $\per{p}$ and $\per{q}$  the
wedge product induces the Poincar\'e Duality isomorphism:
$
\lau{\IH}{*}{\per{p},\sbat}{X}\cong \lau{\IH}{\dim X - *}{\per{q},\sbat}{X}.
$
This gives the  $\R(e)$-isomorphism:
$
\lau{\IL}{}{\per{p}}{X}\cong \lau{\IL}{}{\per{q}}{X}
$
(cf. \cite{Bry}).
It preserves (resp. inverts) the superalgebra structure when $\dim X$ is
even (resp. odd).

\zati The equivariant Gysin sequence  can be written in
the following way
$$
\cdots \to
\lau{\IH}{*}{\per{p}}{B} \otimes \Lambda \eut
\stackrel{\ib{\pi'}{\per{p}}}{\TO}
\lau{\IH}{*}{\per{p},\sbat}{X}
\stackrel{\ib{\oint'}{\per{p}}}{\TO}
\hiru{H}{*}{\lau{\mathcal{G}}{\cdot}{\per{p}}{B}} \otimes \Lambda \eut
\stackrel{\ib{\eub}{\per{p}} \otimes 1+ \ib{I}{\per{p}} \otimes \eut
}{-\!\!\!-\!\!\!-\!\!\!-\!\!\!-\!\!-\!\!\!\TO}
\lau{\IH}{*}{\per{p}}{B} \otimes \Lambda \eut
\to \cdots,
$$
which is a long exact sequence in the category of $\Lambda \eut$-modules.
Since  localization is an exact functor, we get the {\em localized
Gysin sequence}
$$
\cdots \to
\lau{\IH}{*}{\per{p}}{B} \otimes \R(\eut )
\stackrel{\ib{\pi''}{\per{p}} }{\TO}
\lau{\IL}{\star}{\per{p},\sbat}{X}
\stackrel{\ib{\oint''}{\per{p}}}{\TO}
\hiru{H}{*}{\lau{\mathcal{G}}{\cdot}{\per{p}}{B}} \otimes \R(\eut )
\stackrel{\ib{\eub}{\per{p}} \otimes 1+ \ib{I}{\per{p}} \otimes \eut }{-\!\!\!-\!\!\!-\!\!\!-\!\!\!-\!\!-\!\!\!\TO}
\lau{\IH}{*}{\per{p}}{B} \otimes \R(\eut )
\to \cdots,
$$
where
$\ib{\oint''}{\per{p}} ([c ]\ib{\otimes}{\Lambda\euts} R)= \ib{\oint'}{\per{p}} [c] \ib{\otimes}{\Lambda\euts} R$. Thus, we get $\lau{\IL}{\star}{\per{p},\sbat}{X} $ in terms of basic data.

\bigskip

The following result  relates the Euler class
with the  localization of the equivariant intersection cohomology of $X$. It is obtained
straightforwardly from Proposition  \ref{aritz}.

\bp
\label{bi}
Let $X_{1}$, $X_{2}$ be two connected normal unfolded pseudomanifolds.
Consider two modelled actions $\Phi_1 \colon \sbat \times X_1 \to X_1$
and $\Phi_2 \colon \sbat \times X _2\to X_2$.  Let us suppose that
there exists an unfolded isomorphism $f \colon B_{1} \to B_{2}$
between the associated orbit spaces.  Then, the first following
statement implies the second one:

\Zati The isomorphism $f$ is optimal and the Euler classes $e_{1}$
and $e_{2}$ are $f$-related.

\zati There exists a  perverse super\-algebra iso\-mor\-phism
$\mbox{\boldmath $K$} \colon \lau{\mbox{\boldmath $\IL$}}{}{\sbat}{X_{2}}\to
\lau{\mbox{\boldmath $\IL$}}{}{\sbat}{X_{1}}$
verifying
$
\mbox{\boldmath $K$} \rondp \mbox{\boldmath $\pi''_{2}$} = \mbox{\boldmath
$\pi''_{1}$} \rondp (\mbox{\boldmath $f$} \otimes 1 )$.

\ep

The reciprocal to this Theorem does not hold:
just consider the Hopf action on $\S^3$ and the action by multiplication
on the second factor of $\S^2\times\S^1$. The Euler classes are different,
but as the actions are free, both localizations vanish.

\section{\bf Examples}
We illustrate the  results of this work with some particular modelled actions $\Phi
\colon \sbat \times X \to X$. We present: (a) the Gysin and co-Gysin terms, (b) the
equivariant intersection cohomology and (c) the localization of the
equivariant intersection cohomology.

\prg{\bf The pseudomanifold $X$ is a manifold}. Consider the case where $\per{0} \leq \per{p} \leq \per{t}$\footnote{In this range the intersection cohomology of $X$ coincides with its cohomology (see for example \cite{S2}).}.

\Zati $\hiru{H}{*}{\lau{\mathcal{G}}{\cdot}{\per{p} }{B}}=
\lau{\IH}{*}{\per{p}- \per{e}}{B}$ (cf. \cite[sec. 6.4]{Pa}) and $
\lau{H}{*}{}{\lau{\mathcal{K}}{\cdot}{\per{p}}{B}}
=
\lau{H}{*}{\frac{\per{p}}{\per{p} -
\per{e}}}
{B}
$
(cf. \cite{K}). In particular, $\hiru{H}{*}{\lau{\mathcal{G}}{\cdot}{\per{0} }{B}}=
\hiru{H}{*}{B,F}$ and $
\lau{H}{*}{}{\lau{\mathcal{K}}{\cdot}{\per{0}}{B}}
=
\lau{H}{*}{}
{F}$.

\zati $\lau{\IH}{i}{\per{p},\sbat}{X}= \lau{H}{*}{\sbat}{X}$, $\tres{E}{i,0}{\per{p},2} =
\lau{\IH}{i}{\per{p}}{B}$ and
$\displaystyle\tres{E}{i,2j}{\per{p},2} =
\prod_{S \in \mathcal{S}_X}
\hiru{H}{i-2\left[ \frac{\per{p}(S)}{2}\right]}{S}
\otimes \eut^j,
$
for $j>0$.

\zati $\lau{\IL}{\star}{\per{p},\sbat}{X} = \hiru{H}{}{F}\otimes \R(e)
 = \lau{L}{*}{\sbat}{X}.$

\prg{\bf Free action}. These actions are characterized by the condition  $\per{e} = \per{x}=\per{0}$.

\Zati $\hiru{H}{*}{\lau{\mathcal{G}}{\cdot}{\per{p}}{B}}=
\lau{\IH}{*}{\per{p}}{B}$  and $
\lau{H}{*}{}{\lau{\mathcal{K}}{\cdot}{\per{p}}{B}}
=
0$.

\zati The basic spectral sequence  degenerates at the
second term and we have $\lau{\IH}{*}{\per{p},\sbat}{X} = \lau{\IH}{*}{\per{p}}{B}$.
The $\Lambda \eut$-module structure is given by $\eut \cdot b=
e\wedge b$. The perverse super structure comes from that of $B$.

\zati $\lau{\IL}{\star}{\per{p},\sbat}{X} = 0.$
We observe that neither the Euler class $e$ nor the Euler perversity $\per{e}$
are determined by the localization of the intersection cohomology.

\prg{\bf Action without  perverse strata.} These actions are characterized by the condition  $\per{e} = \per{x}$.

\Zati We have
 $\lau{\mathcal{G}}{*}{\per{p}}{B}=
\lau{\Omega}{*}{\per{p} - \per{x}}{B}$ and $
\lau{\mathcal{K}}{*}{\per{p}}{B}
=
\lau{\Om}{*}{\frac{\per{p}}{\per{p} -
\per{x}}}
{B}
$.

\zati
$
\tres{E}{i,2j}{\per{p},2}  =
\left\{
\begin{array}{lc}
\lau{H}{i}{\per{p}}{B} & \hbox{if } j=0 \\[,3cm]
\lau{\IH}{i}{\frac{\per{p}}{\per{p} -\per{x}}}{B}\otimes \eut^j& \hbox{if } j>0
\end{array}
\right.
$.

\zati The perverse  super structure comes from that of $B$.

\prg{\bf The Euler class $e$ is zero}. In particular, all the fixed
strata are non-perverse. We
have

\Zati $\hiru{H}{*}{\lau{\mathcal{G}}{\cdot}{\per{p}}{B}}=
\lau{\IH}{*}{\per{p} - \per{x}}{B}$  and $
\lau{H}{*}{}{\lau{\mathcal{K}}{\cdot}{\per{p}}{B}}
=
\lau{H}{*}{\frac{\per{p}}{\per{p} -
\per{x}}}
{B}
$.

\zati
$\lau{\IH}{*}{\per{p}}{X} =
\lau{\IH}{*}{\per{p}}{B}
\oplus
\lau{\IH}{*-1}{\per{p} - \per{x}}{B}.$
 The basic spectral sequence degenerates at the
second term and we have $\lau{\IH}{*}{\per{p},\sbat}{X} =
\lau{\IH}{*}{\per{p}}{B} \oplus
\left\{ \lau{\IH}{*}{\frac{\per{p}}{\per{p} -
\per{x}}}
{B}  \otimes \bi{\Lambda}{>0}
\eut
\right\}.$
The $\Lambda \eut$-module struture is given by
$
\eut \cdot (b_0,\overline{b_1} \otimes \eut^n) =
(0,\overline{b_0} \otimes \eut + \overline{b_1} \otimes \eut^{n+1}) .$

\zati  $\lau{\IL}{\star}{\per{p},\sbat}{X} = \lau{\IH}{*}{\frac{\per{p}}{\per{p} -
\per{x}}}{B} \otimes \R(\eut).$ The perverse super structure comes from that
of $B$.

\prg{\bf Local calculation}. Consider a chart $(U,\phii)$ of a fixed point $x$
lying on a stratum $S$. The open subset $U$ is $\sbat$-invariant and
describes the local geometry near $x$. It can be equivariantly retracted by isomorphisms
to $c\ib{L}{S}$,  endowed with the action $\Phi_{\ib{L}{S}}$. So, it is enough to consider the case $U = c\ib{L}{S}$. We have

\Zati
$
\hiru{H}{i}{\lau{\mathcal{G}}{\cdot}{\per{p}}{U/\sbat}}=
\left\{
\begin{array}{cl}
\hiru{H}{i}{\lau{\mathcal G}{\cdot}{\per{p}}{\ib{L}{S}/\sbat}}
& \hbox{if } i \leq
m-2
\\[,3cm]
\Ker
\left\{ \ib{\eub}{\per{p}} \colon
\hiru{H}{i}{\lau{\mathcal G}{\cdot}{\per{p}}{\ib{L}{S}/\sbat}}
\TO
\lau{\IH}{i+2}{\per{p}}{\ib{L}{S}/\sbat} \right\}& \hbox{if } i =
m-1
\\[,3cm]
0 &  \hbox{if } i \geq
m.
\end{array}
\right.
$

\nt (cf. \cite[sec. 7.2]{Pa}).

\zati The computation of $\lau{\IH}{*}{\per{p},\sbat}{U}$ is achieved through the following step-by-step procedure.
Let   $\per{q}$ the perversity defined by: $\per{q} = \per{p}$ on $U \backslash S$ and $\per{q}(U \cap S) =
 \per{p}(U \cap S)-1 = m-1$. We have
 $
 \lau{\IH}{i}{\per{p},\sbat}{U} =\lau{\IH}{i}{\per{q},\sbat}{U}
 $, for $i \ne m, m+1$,
 and the exact sequence
$$
 0 \to  \lau{\IH}{m}{\per{q},\sbat}{U} \to  \lau{\IH}{m}{\per{p},\sbat}{U} \to
 \lau{\IH}{m}{\per{p}}{\ib{L}{S}} \otimes \Lambda \eut
  \to \lau{\IH}{m+1}{\per{q},\sbat}{U} \to  \lau{\IH}{m+1}{\per{p},\sbat}{U} \to 0
 $$
For example, when
the action $\Phi_{\ib{L}{S}}$ is free we get that
$$
\lau{\IH}{*}{\per{p},\sbat}{U} =
 \lau{\IH}{\leq m-1}{\per{p}}{\ib{L}{S}/\sbat}
\oplus
\left\{\lau{\IH}{m}{\per{p}}{\ib{L}{S}/\sbat} \otimes
 \Lambda \eut
 \right\}
\oplus
\left\{
\fracc{\lau{\IH}{m-1}{\per{p}}{\ib{L}{S}/\sbat} }{\Ker
\left\{\ib{\eub}{\per{p}} \colon \lau{\IH}{m-1}{\per{p}}{\ib{L}{S}/\sbat} \to
\lau{\IH}{m+1}{\per{p}}{\ib{L}{S}/\sbat} \right\}}\otimes
 \bi{\Lambda}{>0} \eut
\right\} .
$$
 The $\Lambda \eut$-product is induced by
$
\eut \cdot (b_1,\overline{b_2} \otimes \eut^n,b_3 \otimes \eut^n)=
\left\{
\begin{array}{ll}
(b_1 \wedge e, \overline{b_2} \otimes \eut^{n+1} , b_3\otimes \eut^{n+1})
& \hbox{if } |b_1| \leq m-3\\
(0, \overline{b_1 \wedge \eut} \otimes 1 + \overline{b_2} \otimes \eut^{n+1} , b_3\otimes
\eut^{n+1})
& \hbox{if } |b_1|  = m-2\\
(0,  \overline{b_2} \otimes \eut^{n+1} , b_1 \otimes \eut +b_3\otimes
\eut^{n+1})
& \hbox{if } |b_1|  = m-1.
\end{array}
\right.
$

\zati We have a long exact sequence
$$
\cdots
\to
\lau{\IL}{\star}{\per{q},\sbat}{U}
\to
\lau{\IL}{\star}{\per{p},\sbat}{U}
\to
\lau{\IH}{m}{\per{p}}{\ib{L}{S}} \otimes \R(e)
\to
\lau{\IL}{\star}{\per{q},\sbat}{U}
\to
\lau{\IL}{\star}{\per{p},\sbat}{U}
\to
\cdots.
$$
When  $\Phi_{\ib{L}{S}}$ is free then
\begin{equation}
\label{amaia}
\lau{\IL}{\star}{\per{p},\sbat}{U} =
\left\{
\fracc{\lau{\IH}{m-1}{\per{p}}{\ib{L}{S}/\sbat} }{\Ker \left\{\ib{\eub}{\per{p}} \colon \lau{\IH}{m-1}{\per{p}}{\ib{L}{S}/\sbat} \to
\lau{\IH}{m+1}{\per{p}}{\ib{L}{S}/\sbat} \right\}}
\oplus
\lau{\IH}{m}{\per{p}}{\ib{L}{S}/\sbat} \right\}\otimes
\R(\eut) .
\end{equation}
The perverse super structure comes from that of $\ib{L}{S}/\sbat$.

\end{document}